# The Repeated Solicitation Model

**R. W. R. Darling***

## ABSTRACT

* P. O. Box 535, Annapolis Junction, Maryland 20701-0535, USA

**e-mail**: rwrd@afterlife.ncsc.mil

This paper presents a probabilistic analysis of what we call the "repeated solicitation model". To give a specific context, suppose B is a direct marketing company with a list of $S_0$ sales prospects. At epoch 1, B sends a solicitation to every prospect on the list, and elicits $X_1$ replies. The company deletes the respondents from the list, and at epoch 2 sends a solicitation to the other prospects, of whom $X_2$ respond, and so on. This continues until an epoch $n$ such that $X_n = 0$, which we call epoch $T$, and then B makes no further solicitations. We seek

- The probability distribution of $T$;
- The distribution of the total number of respondents, $X_1 + \ldots + X_{T-1}$;
- The expected total number of solicitations.

All three quantities are explicitly computed, assuming that (a) prospects' response times are independent, and (b) $S_0$ is Poisson distributed.

# 1.0 Definition and Principal Results

## 1.1 The Repeated Solicitation Model

Consider a collection of independent random variables $\{U_1, \ldots, U_r\}$ with values in $\{1, 2, 3, \ldots, \infty\}$ (the events $\{U_s = \infty\}$ are possible), all with the same distribution, specified by

$$\pi_n \equiv P(U_1 = n)\,,\; n \in \{1, 2, \ldots, \infty\}\,. \tag{1}$$

Think of $U_s$ as the first time (which may be never) at which a "client" $s$ would respond, if she were to receive an infinite sequence of "solicitations" at times 1, 2, … . We imagine a "solicitor", B, who solicits each of his clients at the times 1, 2, … , except that no more solicitations are sent to a client after that client has responded. For $n \in \{1, 2, \ldots, \infty\}$, define

$$X_n \equiv \sum_{s=1}^{r} 1_{\{U_s = n\}}\,. \tag{2}$$

Interpret $X_n$ as the number of responses which B receives at time $n$ from the pool of clients, while $X_\infty$ is the number who never respond. This also has an urn model interpretation - see **[2]** for some related constructions. However at the first time $n$ in which no responses are received, B goes into a state of "despair" and abandons the solicitation process. Three random variables are of interest:

$$T \equiv \min\{n \ge 1 : X_n = 0\}\,, \tag{3}$$

meaning the time at which B goes into despair;

$$Y \equiv X_1 + \ldots + X_{T-1} = \sum_{s=1}^{r} 1_{\{U_s < T\}}\,, \tag{4}$$

meaning the number of responses received by B before going into despair; and

$$M \equiv \sum_{s=1}^{r} \left(U_s 1_{\{U_s < T\}} + T 1_{\{U_s > T\}}\right)\,, \tag{5}$$

meaning the number of solicitations made by B before going into despair.

## 1.2 Short Cut: the Poisson-Multinomial Relationship

First, observe that the joint law of $X_1, X_2, \ldots, X_\infty$ is given by:

$$(X_1, X_2, \ldots, X_\infty) \sim \text{Multinomial}(r; \pi_1, \pi_2, \ldots, \pi_\infty) . \tag{6}$$

Instead of taking a fixed client base, we will assume that the initial number of clients $S_0 \sim \text{Poisson}(\nu)$, for some parameter $\nu$. This is natural in the marketing context, because the set of clients may well be selected at random from some larger data base. The author thanks Mark Jacobson for pointing out the following basic fact, which allows us to present simple proofs for most of the main results.

### 1.2.1 Lemma

*Under the prior distribution* $S_0 \sim \text{Poisson}(\nu)$, $X_1, X_2, \ldots, X_\infty$ *become independent Poisson random variables, where* $X_n$ *has parameter* $\pi_n \nu$.

**Proof:** Take any sequence of non-negative integers $(x_1, x_2, \ldots, x_\infty)$, all but finitely many of which are zero. Multiply the multinomial probability $P(X_1 = x_1, X_2 = x_2, \ldots, X_\infty = x_\infty | S_0 = r)$ by the Poisson probability $P(S_0 = r)$. Sum over $r$ such that $r$ is the sum of the $\{x_i\}$ (there is only one!). The result factors as a product of Poisson probabilities, as desired. o

## 1.3 Main Results

For convenience, introduce parameters $\{\lambda_n\}$ given by:

$$\lambda_n \equiv 1 - e^{-\pi_n \nu} = P(X_n > 0) . \tag{7}$$

We know $P(T < \infty) = 1$, since $P(T < S_0 + 1) = 1$. Lemma [1.2.1], allows us to write down the distribution of $T$ when $S_0 \sim \text{Poisson}(\nu)$, namely

$$P(T = n) = P(X_1 > 0) \ldots P(X_{n-1} > 0) P(X_n = 0) = \lambda_1 \ldots \lambda_{n-1}(1 - \lambda_n) . \tag{8}$$

The identity $\mathrm{E}[T] = \sum_{n \ge 1} P(T \ge n)$ implies

$$\mathrm{E}[T] = 1 + \sum_{n \ge 1} \lambda_1 \ldots \lambda_n . \tag{9}$$

Here is our main result, which depends crucially on the choice of prior distribution for $S_0$.

### 1.3.1 Theorem

*If* $S_0 \sim \text{Poisson}(\nu)$*, then conditional on* $T = k \geq 2$*,* $Y$ *may be expressed as a sum* $X_1' + \ldots + X_{k-1}'$ *of independent, strictly positive random variables, with*

$$P(X_n' = m) = P(X_n = m)/\lambda_n,\ m \geq 1. \tag{10}$$

*If F is the distribution function of* $U_1$*, i. e.* $F(k) \equiv \pi_1 + \ldots + \pi_k$*, then*

$$\mathrm{E}[Y] = \nu\mathrm{E}[F(T)]; \tag{11}$$

$$\mathrm{E}[M] = \nu\mathrm{E}[H(T)], \tag{12}$$

*where* $H(k) \equiv \mathrm{E}[U_1 \wedge k]$*. Moreover*

$$\mathrm{Var}(Y) = \nu\mathrm{E}[J(T)] + \nu^2\mathrm{Var}(R(T)), \tag{13}$$

*where*

$$J(k) \equiv \sum_{n=1}^{k} \pi_n \left\{1 - \frac{\pi_n \nu(1-\lambda_n)}{\lambda_n}\right\}, \tag{14}$$

$$R(1) \equiv 0,\ R(k) \equiv \frac{\pi_1}{\lambda_1} + \ldots + \frac{\pi_{k-1}}{\lambda_{k-1}},\ k \geq 2. \tag{15}$$

**Proof: Part I.** Suppose $S_0 \sim \text{Poisson}(\nu)$. According to the independence assertion of Lemma [1.2.1], we see that $X_n$ is independent of $\{T \leq n-1\}$, and hence of its complement $\{T \geq n\}$. Arguing as in the proof of Wald's equation (**[1]**, p. 396), we see that

$$\mathrm{E}Y = \mathrm{E}[X_1 + \ldots + X_T] = \sum_{n \geq 1} \mathrm{E}[X_n 1_{\{T \geq n\}}] = \sum_{n \geq 1} \mathrm{E}[X_n] P(T \geq n)$$

$$= \sum_{n \geq 1} \pi_n \nu \sum_{k \geq n} P(T = k) = \nu \sum_{k \geq 1} P(T = k) \sum_{n=1}^{k} \pi_n,$$

which establishes [11].

**Part II**. To establish the representation of *Y*, note that, for $n < k$ and $m \geq 1$, the independence assertions of Lemma [1.2.1] imply

$$P(X_n = m | T = k) = P(X_n = m | X_1 > 0, \ldots, X_{k-1} > 0, X_k = 0)$$

$$= P(X_n = m | X_n > 0) = \frac{P(X_n = m)}{\lambda_n}. \quad \textbf{(16)}$$

It follows that, conditional on $T = k$, we may write

$$Y = X_1' + \ldots + X_{k-1}' \quad \textbf{(17)}$$

where $X_1', X_2', \ldots$ are conditionally independent with the law [10]. (The last line can be used to give another proof of [11], along the lines of [19].)

**Part III**. Next we will prove [12]. We may write *M* in the form

$$M = \sum_{n=1}^{\infty} \{n1_{\{T \geq n\}} + T1_{\{T < n\}}\} X_n,$$

where $n = \infty$ is included in the summation. As noted above, $X_n$ is independent of $T1_{\{T \leq n-1\}}$, and of $\{T \geq n\}$. Arguing as in Step I,

$$\mathrm{E}M = \sum_{n=1}^{\infty} \{nP(T \geq n) + \mathrm{E}[T1_{\{T < n\}}]\}\mathrm{E}[X_n]$$

$$= \nu \sum_{n \geq 1} \pi_n \sum_{k \geq 1} (k \wedge n) P(T = k)$$

$$= \nu \sum_{k \geq 1} P(T = k) \sum_{n \geq 1} \pi_n (k \wedge n),$$

which verifies [12].

**Part IV**. Finally we establish [13], on the basis of the formula

$$\mathrm{Var}(Y) = \mathrm{Var}(\mathrm{E}[Y|T]) + \mathrm{E}[\mathrm{Var}(Y|T)]. \tag{18}$$

As for the first summand of [18], [17] implies that

$$\mathrm{E}[Y|T = k] = \nu R(k) \equiv \nu\left\{\frac{\pi_1}{\lambda_1} + \ldots + \frac{\pi_{k-1}}{\lambda_{k-1}}\right\}.$$

and so

$$\mathrm{Var}(\mathrm{E}[Y|T]) = \nu^2 \mathrm{Var}(R(T)).$$

By the conditional independence of the summands in [17],

$$\mathrm{Var}(Y|T = k) = \sum_{n=1}^{k-1} \mathrm{Var}(X_n').$$

For brevity, write $b_n \equiv \mathrm{Var}(X_n')$. Since $P(T > n) = \lambda_n P(T > n-1)$ by [8], Fubini's Theorem shows that

$$\sum_{k \geq 2} P(T = k) \sum_{n=1}^{k-1} b_n = \sum_{n \geq 1} b_n \sum_{k \geq n+1} P(T = k) = \sum_{n \geq 1} b_n \lambda_n P(T > n-1)$$

$$= \sum_{n \geq 1} b_n \lambda_n \sum_{k \geq n} P(T = k) = \sum_{k \geq 1} P(T = k) \sum_{n=1}^{k} b_n \lambda_n, \tag{19}$$

or in other words,

$$\mathrm{E}[\mathrm{Var}(Y|T)] = \sum_{k \geq 2} P(T = k)\mathrm{Var}(Y|T = k) = \sum_{k \geq 1} P(T = k) \sum_{n=1}^{k} \mathrm{Var}(X_n') \lambda_n.$$

We compute $\mathrm{Var}(X_n')$, from $\mathrm{E}[X_n] = \pi_n \nu = \mathrm{Var}(X_n)$, giving

$$\sum_{n=1}^{k} \mathrm{Var}(X_n') \lambda_n = \nu \sum_{n=1}^{k} \pi_n\left\{1 - \frac{\pi_n \nu}{\lambda_n}(1 - \lambda_n)\right\} = \nu J(k).$$

and the result [13] follows from [18]. o

**1.3.2 Remark**

Formula [8] makes it practical to compute $\mathrm{E}[f(T)]$ for any function $f$ for which the expectation is finite, and indeed we shall give examples later of the computation of $\mathrm{E}[Y]$ when the law of $U_1$ is geometric.

## 2.0 A Basic Technical Tool for Model Reduction

This section develops a technical tool for deeper analysis of the Repeated Solicitation Model. The set of integers $\{1, 2, 3, \ldots\}$ will be denoted $Z_+$. We shall allow the cardinality $S_0$ of the set of random variables $\{U_s\}$ to vary; we are no longer assuming a Poisson prior for $S_0$.

**2.1 Lemma**

*For any* $g: Z_+ \to [0, \infty)$,

$$\mathrm{E}\left[\sum_{n=1}^{T-1} g(n) X_n \middle| S_0 = r\right] = r\mathrm{E}[f(T) | S_0 = r-1] \tag{20}$$

*(both sides may be infinite), where*

$$f(k) \equiv \sum_{n=1}^{k} g(n)\pi_n = \mathrm{E}[g(U_1) 1_{\{U_1 \le k\}}]. \tag{21}$$

**Proof:** If $S_0 = r \ge 1$, which is omitted from the notation, the left side of [20] can be written as

$$\mathrm{E}\left[\sum_{n=1}^{T-1} g(n) \sum_{s=1}^{r} 1_{\{U_s = n\}}\right] = r \sum_{n \ge 1} g(n) P(T > U_1 = n),$$

since the events $\{U_s < T\}$, for $s \in \{1, 2, \ldots, r\}$, are exchangeable. Now

$$\sum_{n \ge 1} g(n) P(T > U_1 = n) = \sum_{n \ge 1} g(n) P\left(\{U_1 = n\} \cap \bigcap_{m=1}^{n-1} \left(\bigcup_{s=2}^{r} \{U_s = m\}\right)\right)$$

$$= \sum_{n \ge 1} g(n) P(\{U_1 = n\} \cap \{\hat{T} \ge n\}),$$

where $\hat{T}$ is the analog to $T$, with respect to $\{U_2, \ldots, U_r\}$; in other words

$$\hat{T} \equiv \min\left\{n \geq 1: \sum_{s=2}^{r} 1_{\{U_s = n\}} = 0\right\}. \quad \textbf{(22)}$$

Now $P(\{U_1 = n\} \cap \{\hat{T} \geq n\}) = P(U_1 = n)P(\hat{T} \geq n)$, since $\hat{T}$ is independent of $U_1$, so we obtain

$$\sum_{n \geq 1} g(n)P(U_1 = n) \sum_{k=n}^{\infty} P(T = k | S_0 = r-1)$$

$$= \sum_{k \geq 1} P(T = k | S_0 = r-1) \sum_{n=1}^{k} g(n)\pi_n ,$$

by Fubini's Theorem, on reversing the order of summation. **o**

The following result leads to an alternative proof of [11] and [12] above.

**2.2 Proposition**

*If* $S_0 \sim \text{Poisson}(\nu)$*, then for* $g: Z_+ \to [0, \infty)$,

$$\mathrm{E}\left[\sum_{n=1}^{T-1} g(n)X_n\right] = \nu \mathrm{E}[f(T)], \quad \textbf{(23)}$$

*where f and g are related as in [21].*

**Proof:** Since $S_0 \sim \text{Poisson}(\nu)$, Lemma [2.1] gives:

$$\mathrm{E}\left[\sum_{n=1}^{T-1} g(n)X_n\right] = \sum_{r \geq 1} \frac{e^{-\nu}\nu^r}{r!} r\mathrm{E}[f(T)|S_0 = r-1]$$

$$= \nu \sum_{m \geq 0} \frac{e^{-\nu}\nu^m}{m!} \mathrm{E}[f(T)|S_0 = m]$$

which gives the result. **o**

# 3.0 The Markov Model

This section develops another perspective on the Repeated Solicitation Model.

## 3.1 Model Parameters

The probability distribution of $U_1$ can be described in terms of the parameters

$$q_{n-1} = P(U_1 > n | U_1 > n-1),\ n = 1, 2, \ldots, \quad (24)$$

and $p_m \equiv 1 - q_m$. In other words, $\pi_1 \equiv p_0 = P(U_1 = 1)$, and

$$\pi_{n+1} \equiv p_n q_{n-1} \ldots q_0,\ n \geq 1. \quad (25)$$

## 3.2 Number of Non-responding Clients as a Markov Chain

Let us use the recursive construction $S_0 = r$ and

$$S_n \equiv S_{n-1} - X_n,\ n = 1, 2, \ldots,$$

to define a stochastic process $\{S_n\}$, and sigma-fields $\Im_n \equiv \sigma\{S_0, \ldots, S_n\}$ for $n = 0, 1, \ldots$. Observe that there are $S_{n-1}$ clients left in the pool when B solicits in epoch $n$, each of which has a probability $p_{n-1}$ of responding at this epoch; they are conditionally independent given $\Im_{n-1}$, and therefore

$$X_n \sim \text{Binomial}(S_{n-1}, p_{n-1}).$$

It follows that, conditional on $\Im_{n-1}$, $S_n \sim \text{Binomial}(S_{n-1}, q_{n-1})$, so $\{S_n, \Im_n\}$ is a non-homogeneous Markov chain, and $T$ is an $\{\Im_n\}$-stopping-time.

## 3.3 A Sequence of Stopping-Times

It will be useful to define a collection $T \equiv T^{(0)}, T^{(1)}, \ldots$ of $\{\Im_n\}$-stopping-times:

$$T^{(n)} \equiv \min\{k > n : X_k = 0\}.$$

We are interested in the quantities

$$\phi_r^{(n)}(f) \equiv \mathrm{E}[f(T^{(n)}) | S_n = r], \ n \geq 1, \qquad \textbf{(26)}$$

for $f: Z_+ \to [0, \infty)$, which admit the following recursive formula.

**3.3.1 Lemma**

*For any* $f: Z_+ \to [0, \infty)$, $\phi_0^{(n)}(f) = f(n+1)$, *and*

$$\phi_r^{(n)}(f) = f(n+1)q_n^r + \sum_{k=0}^{r-1} \binom{r}{k} p_n^{r-k} q_n^k \phi_k^{(n+1)}(f), \ r \geq 1. \qquad \textbf{(27)}$$

**Proof:** If $S_n = 0$, then $X_{n+1} = 0$ and $T^{(n)} = n+1$, giving $\phi_0^{(n)} = f(n+1)$. For $r \geq 1$,

$$\phi_r^{(n)}(f) = \mathrm{E}\Big[f(T^{(n)}) 1_{\{X_{n+1} = 0\}} \Big| S_n = r\Big] + \sum_{k=0}^{r-1} \mathrm{E}\Big[f(T^{(n)}) 1_{\{X_{n+1} = r-k\}} \Big| S_n = r\Big]$$

$$= f(n+1)q_n^r + \sum_{k=0}^{r-1} \mathrm{E}[f(T^{(n+1)}) | S_{n+1} = k] \binom{r}{k} p_n^{r-k} q_n^k,$$

using the Markov property of $\{S_n\}$, and the result follows. **o**

For $f: Z_+ \to [0, \infty)$ as in [26], define

$$\Phi_\nu^{(n)}(f) \equiv \sum_{r \geq 0} \frac{e^{-\nu} \nu^r}{r!} \phi_r^{(n)}(f), \qquad \textbf{(28)}$$

so that, for example, when $S_0 \sim \mathrm{Poisson}(\nu)$, then

$$\Phi_\nu^{(0)}(f) = \mathrm{E}[f(T)].$$

The special role of the $\mathrm{Poisson}(\nu)$ prior for $S_0$ is clarified somewhat by the following Proposition.

**3.4 Proposition**

*If* $S_0 \sim \mathrm{Poisson}(\nu)$, *then for* $f: Z_+ \to [0, \infty)$ *the* $\Phi_\nu^{(n)}(f)$ *satisfy the recursion*

$$\Phi_\nu^{(n)}(f) = e^{-p_n\nu} f(n+1) + (1 - e^{-p_n\nu})\Phi_{q_n\nu}^{(n+1)}(f). \tag{29}$$

**Proof:** Let us substitute [27] into the right side of [28]. We obtain

$$e^\nu \Phi_\nu^{(n)}(f) = \phi_0^{(n)}(f) + \sum_{r \ge 1} \frac{\nu^r}{r!}\left\{ f(n+1) q_n^r + \sum_{k=0}^{r-1} \binom{r}{k} p_n^{r-k} q_n^k \phi_k^{(n+1)}(f) \right\}$$

$$= e^{q_n\nu} f(n+1) + \sum_{k \ge 0} \frac{(q_n\nu)^k}{k!} \phi_k^{(n+1)}(f) \sum_{r \ge k+1} \frac{(p_n\nu)^{r-k}}{(r-k)!}$$

$$= e^{q_n\nu} f(n+1) + e^{q_n\nu} \Phi_{q_n\nu}^{(n+1)}(f)(e^{p_n\nu} - 1)$$

and the result follows. o

This leads to a computable series expansion for $\mathrm{E}[f(T)]$, giving an alternative derivation of [8].

### 3.5 Proposition

*If* $S_0 \sim \mathrm{Poisson}(\nu)$, *then* for *any* $f: Z_+ \to [0, \infty)$,

$$\mathrm{E}[f(T)] = (1 - \lambda_1) f(1) + \sum_{n \ge 2} \lambda_1 \ldots \lambda_{n-1}(1 - \lambda_n) f(n) \tag{30}$$

*where* $\{\lambda_n\}$ *are as in [7].*

**Proof:** Applying [29] recursively, we find that

$$\Phi_\nu^{(0)}(f) = e^{-p_0\nu} f(1) + (1 - e^{-p_0\nu})\left\{ e^{-p_1 q_0 \nu} f(2) + (1 - e^{-p_1 q_0 \nu})\Phi_{q_1 q_0 \nu}^{(2)}(f) \right\},$$

$$= e^{-\pi_1\nu} f(1) + (1 - e^{-\pi_1\nu}) e^{-\pi_2\nu} f(2) + (1 - e^{-\pi_1\nu})(1 - e^{-\pi_2\nu}) e^{-\pi_3\nu} f(3) + \ldots$$

In the case where $f(k) = 1_{\{k = m\}}$, the only non-zero term in the series is $\lambda_1 \ldots \lambda_{m-1}(1 - \lambda_m)$, while the left side is $P(T = m)$, as desired. Now [30] follows since $f$ is non-negative. o

## 4.0 Special Case: Geometric Response Times

A simple model in the marketing context might assign to each of the $U_i$ a distribution which is a mixture of: a point mass at 1, reflecting pent-up demand; a point mass at infinity, reflecting total indifference; and a Geometric($p$) distribution on the positive integers. For brevity, we shall only present the results for the pure Geometric case, i.e. where $p_n = p$ and $q_n = q = 1 - p$ are constant in [24], and [1] is replaced by

$$P(U_s > m) = q^m . \quad (31)$$

The Markov chain $\{S_n, \mathfrak{I}_n\}$ is now homogeneous, so the law of $T^{(n)} - n$, given that $S_n = r$, has the same law as $T$ given that $S_0 = r$. A key role is played by the probability generating function

$$G_\nu(z) \equiv \mathrm{E}[z^T], \quad (32)$$

which is given by [30], and satisfies

$$\frac{G_\nu(z)}{z} = e^{-p\nu} + (1 - e^{-p\nu}) G_{q\nu}(z) . \quad (33)$$

### 4.1 Expected Total Yield

#### 4.1.1 Proposition

*If* $U_1 \sim \mathrm{Geometric}(p)$ *and* $S_0 \sim \mathrm{Poisson}(\nu)$*, then* $\mathrm{E}[Y]$ *is explicitly computable from the formula*

$$\mathrm{E}[Y] = \nu\{1 - G_\nu(q)\} . \quad (34)$$

**Proof:** This follows from [11], [31], and [32]. **o**

### 4.2 Expected Total Marketing Effort

In the case of constant $p$, a special relationship holds between $\mathrm{E}[M]$ and $\mathrm{E}[Y]$.

### 4.3 Proposition

*If* $U_1 \sim \text{Geometric}(p)$*, the expected total marketing effort is given by*

$$\mathrm{E}[M|S_0 = r] = \mathrm{E}[Y|S_0 = r]/p\,. \tag{35}$$

**Proof:** Although simpler proofs are possible, we shall give a construction from which further information can also be deduced. Consider the $\{\Im_n\}$-adapted stochastic process $\{V_n, n \ge 0\}$ given by $V_0 \equiv 1$, and

$$V_n \equiv \frac{z^{X_1 + \ldots + X_n}}{(q + pz)^{S_0 + \ldots + S_{n-1}}},\; n \ge 1\,.$$

where $z$ is an arbitrary positive real number. For fixed $z$, each $V_n$ is a bounded random variable, and $\{(V_n, \Im_n), n \ge 0\}$ is a martingale, since

$$\mathrm{E}[z^{X_n}|\Im_{n-1}] = (q + pz)^{S_{n-1}}.$$

When $S_0 = r$, the stopping-time $T$ is bounded above by $r + 1$. The conditions of the Optional Stopping Theorem (**[1]**, p. 464) hold, showing that

$$1 = \mathrm{E}[V_0] = \mathrm{E}[V_T] = \mathrm{E}\left[\frac{z^Y}{(q + pz)^M}\right].$$

using [5]. Differentiate with respect to $z$, and set $z \equiv 1$, to obtain [35]. **o**

#### 4.3.1 Corollary

*If* $U_1 \sim \text{Geometric}(p)$ *and* $S_0 \sim \text{Poisson}(\nu)$*, then the expected total marketing effort is explicitly computable from the formula*

$$\mathrm{E}[M] = \nu\{1 - G_\nu(q)\}/p\,. \tag{36}$$

**Proof:** Immediate from [34] and [35]. **o**

### 4.4 Formulas Using a Degenerate or Binomial Prior

We continue to suppose $U_1 \sim \text{Geometric}(p)$, but now assume $S_0 \sim \text{Binomial}(s, \theta)$ (which includes the degenerate case where $\theta = 1$). The results of this section do not appear to lead to computable expres-sions; they are included only for the sake of completeness. The counterpart to [27] with $n = 0$ and $f(k) \equiv z^k$ reads

$$g_r(z) \equiv \mathrm{E}[z^T | S_0 = r] = z\left\{q^r + \sum_{k=0}^{r-1} \binom{r}{k} q^k p^{r-k} g_k(z)\right\},\ r \geq 1. \quad \textbf{(37)}$$

Given $\theta \in [0, 1]$, let $G_{s,\theta}(z)$ be the probability generating function of $T$ when $S_0 \sim \mathrm{Binomial}(s, \theta)$, i. e.

$$G_{s,\theta}(z) \equiv \mathrm{E}[\mathrm{E}[z^T | S_0]] = \sum_{k=0}^{s} \binom{s}{k} \theta^k (1-\theta)^{s-k} g_k(z),\ z \geq 0. \quad \textbf{(38)}$$

### 4.5 Lemma

*Taking* $q \equiv 1 - p$, $G_{s,\theta}(z)$ *satisfies the functional equation:*

$$G_{s,\theta}(z) = z\left\{G_{s,\theta q}(z) + (1-\theta p)^s [1 - G_{s,\theta q/(1-\theta p)}(z)]\right\}. \quad \textbf{(39)}$$

*In particular, since* $g_s(z) = G_{s,1}(z)$, *the case* $\theta = 1$ *gives*

$$g_s(z) = \frac{z}{1 + z q^s}\left\{q^s + G_{s,q}(z)\right\}. \quad \textbf{(40)}$$

#### 4.5.1 Remark

As $\theta \to 0$ and $s \to \infty$ such that $\theta s = \nu$, equation [39] becomes equation [33] in the limit.

**Proof:** Combining [37] and [38], and the identity $\binom{s}{r}\binom{r}{k} = \binom{s}{k}\binom{s-k}{r-k}$,

$$\frac{G_{s,\theta}(z)}{z} = \sum_{r=0}^{s} \binom{s}{r} \theta^r (1-\theta)^{s-r} \left\{q^r + \sum_{k=0}^{r-1} \binom{r}{k} q^k p^{r-k} g_k(z)\right\}$$

$$= (1-\theta p)^s + \sum_{k=0}^{s-1} \binom{s}{k} (\theta q)^k g_k(z) \sum_{r=k+1}^{s} \binom{s-k}{r-k} (\theta p)^{r-k} (1-\theta)^{s-r}$$

$$= (1-\theta p)^s + \sum_{k=0}^{s-1} \binom{s}{k} (\theta q)^k g_k(z)\{(1-\theta q)^{s-k} - (1-\theta)^{s-k}\}.$$

Note that, when $k = s$, the factor $\{\ldots\}$ vanishes, so we may write

$$= (1-\theta p)^s + \sum_{k=0}^{s} \binom{s}{k} (\theta q)^k g_k(z)\{(1-\theta q)^{s-k} - (1-\theta)^{s-k}\}$$

$$= (1-\theta p)^s + G_{s,\theta q}(z) - (1-\theta p)^s \sum_{k=0}^{s} \binom{s}{k} \left(\frac{\theta q}{1-\theta p}\right)^k \left(\frac{1-\theta}{1-\theta p}\right)^{s-k} g_k(z)$$

$$= (1-\theta p)^s \{1 - G_{s,\theta q/(1-\theta p)}(z)\} + G_{s,\theta q}(z),$$

as desired. o

**4.5.2 Lemma**

If $U_1 \sim \text{Geometric}(p)$ and $S_0 \sim \text{Binomial}(s, \theta)$, *then*

$$\mathrm{E}[Y] = s\theta\{1 - G_{s-1,\theta}(q)\}.$$

**Proof:** This follows from Lemma [2.1] and [31]. o

## 5.0 Applications to Direct Marketing

We return to the marketing context described in the Abstract. Throughout this section, suppose that a known proportion $\theta$ of those who respond actually purchase the product. For simplicity, we will treat the case where $U_1 \sim \text{Geometric}(p)$ and $S_0 \sim \text{Poisson}(\nu)$. We present here three problems which our model can help to solve, all based on [34] and the formula (from [30])

$$G_\nu(q) = (1-\lambda_1)q + \sum_{n \ge 2} \lambda_1 \ldots \lambda_{n-1}(1-\lambda_n)q^n, \quad \lambda_i \equiv 1 - e^{-pq^{i-1}\nu}. \qquad \textbf{(41)}$$

### 5.1 Selecting the Size of the Prospect Pool

Suppose that $p$ is fixed, and the marketer aims to sell at least $K$ items. The goal is to choose the smallest parameter $\nu$ (i.e. expected size of the num-ber of sales prospects) such that expected sales are at least $K$, i. e. such that, if $E_{\nu,\,p}[Y]$ is the expected total yield [34], then

$$\theta E_{\nu,\,p}[Y] \geq K\ . \tag{42}$$

This can be computed by plotting $E_{\nu,\,p}[Y]$ against $\nu$, using [34] and [41].

### 5.1.1 Example

For $p = 1/512$, we display $E_{\nu,\,p}[Y]$ versus $\nu$, for $0 \leq \nu \leq 3000$. (see graph).

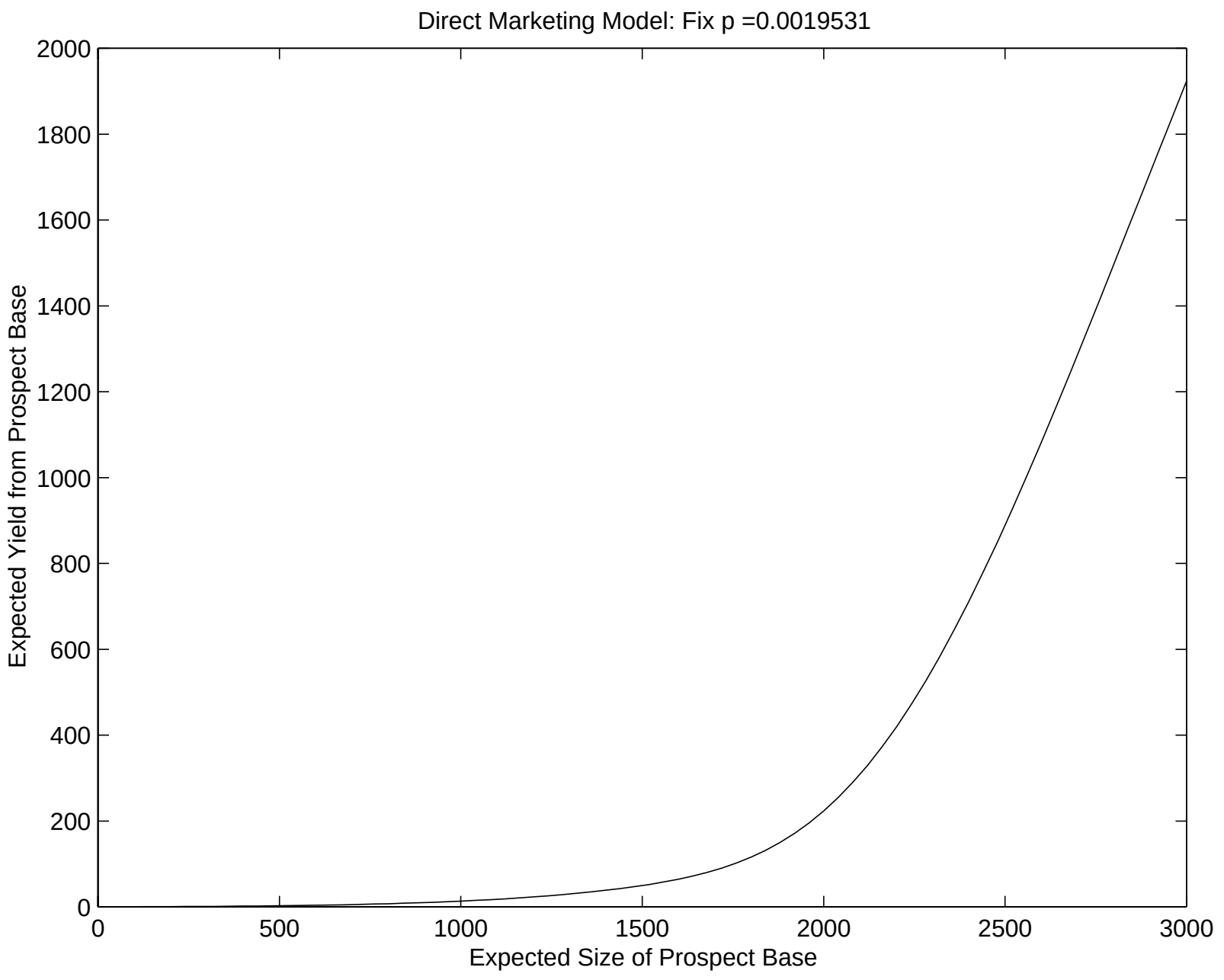

### 5.1.2 How Many Terms to Use in the Expansion

In using the expansion [41] to compute $G_\nu(q)$, we should select a small number $\alpha$, and choose the number of terms $n$ to be so large that

$$P(T > n) \leq \alpha\ .$$

This is easy to do because the law of $T$ is given by [8]. Thus, for $\lambda_j$ as in [41], we stop adding terms as soon as we reach an $n$ such that

$$\lambda_1 \ldots \lambda_n \leq \alpha\ .$$

### 5.2 Selecting Response Probability, for a Fixed Size of the Prospect Pool

Suppose that $\nu$ is fixed, and the marketer aims to sell at least $K$ items. The goal is to choose the smallest response probability $p$ (by making the offer more or less attractive) such that expected sales are at least $K$. This can be computed by plotting $E_{\nu,p}[Y]$ against $p$, using [34] and [41].

#### 5.2.1 Example

Take $\nu = 1,000$; $E_{\nu,p}[Y]$ is shown on the vertical axis, versus $p$ on the horizontal axis (log scale), for $2^{-10} \le p \le 2^{-8}$.

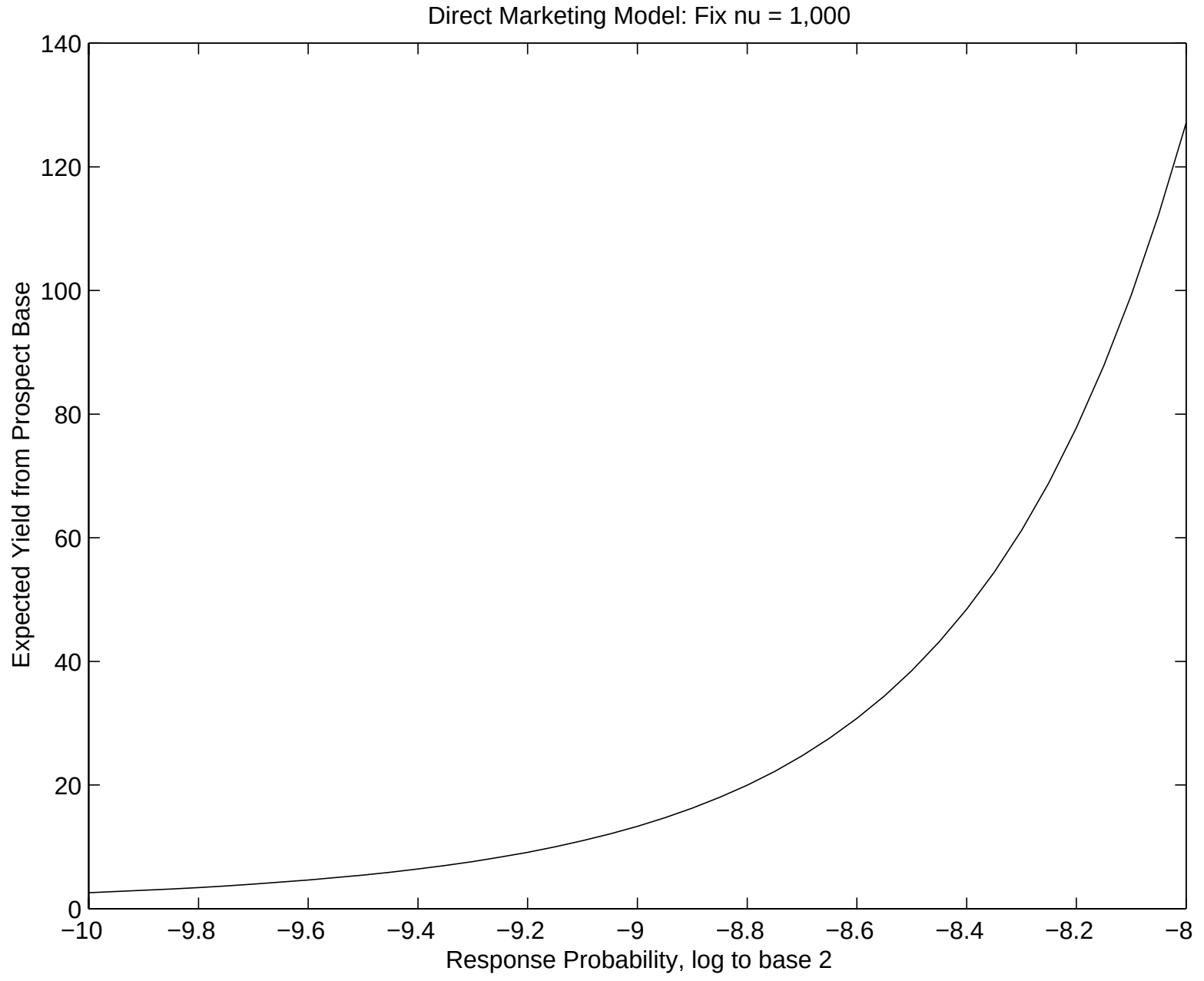

### 5.3 Selecting Two Parameters Simultaneously

Suppose now that there are two variables under the marketer's control, namely the selling price, $w$, and the expected number of prospects, $\nu$. Clearly one could increase sales either by lowering $w$ or by increasing $\nu$, both of which reduce the profit per sale. The parameter $p$, denoting a prospect's probability of responding at each round, is a decreasing function $p(w)$ of $w$. There are three costs faced by the marketer:

- The cost the marketer pays for the product, denoted $c_0$;
- The cost of acquiring each sales prospect, denoted $c_1$;

- The cost of each solicitation, denoted $c_2$.

The expected profit $F(\nu, w)$ is the expected gross revenue from sales, minus the expected costs, namely

$$F(\nu, w) = \theta(w - c_0)\mathrm{E}_{\nu, p(w)}[Y] - c_1\nu - c_2\mathrm{E}_{\nu, p(w)}[M],$$

where $\theta$ is the proportion of respondents who actually purchase the product, and $\mathrm{E}_{\nu, p(w)}[M]$ is the expected total marketing effort [36], which depends on $\nu$ and $p(w)$. For specific choices of parameters $\theta, c_0, c_1, c_2$ and of the function $p(w)$, the formulas given in this paper allow $F(\nu, w)$ to be computed explicitly, so that an optimal pair $(\nu, w)$ can be selected.

## 6.0 Unsolved Problems

Suppose $S_0 = r$, and $r$ is large. It is impractical to compute $P(T \le n)$ by inclusion-exclusion, since there are too many terms. Nevertheless is there a computable formula for $P(T = n)$? Lemma [4.5] gives a functional equation satisfied by the probability generating function (p.g.f.) of $T$, when $U_1$ is geometric; even in this simple case, no computable formula for the p.g.f. results. Alternatively, is there any prior distribution for $S_0$ other than the Poisson which allows $\mathrm{E}[f(T)]$ to be computed using a recursion like that of Proposition [3.5]?

**Acknowledgments:** The author thanks Clark T. Benson for suggesting the problem, and Mark T. Jacobson for illuminating insights and for correcting mathematical errors.